\documentclass[11pt]{amsart}
\usepackage{amssymb, latexsym}
\theoremstyle{plain}
\newtheorem{theorem}{Theorem}

\newtheorem{proposition}{Proposition}

\theoremstyle{remark}

\newtheorem*{Remark 1}{Remark 1}
\newtheorem*{Remark 2}{Remark 2}
\newtheorem*{Remark 3}{Remark 3}
\newtheorem*{Remark 4}{Remark 4}

\numberwithin{equation}{section}

\begin{document}

\title[Diffusions That Eventually Stop Down-Crossing]
 {One-Dimensional Diffusions That Eventually Stop Down-Crossing}

\author{Ross G. Pinsky}
\address{Department of Mathematics\\
Technion---Israel Institute of Technology\\
Haifa, 32000\\ Israel} \email{pinsky@math.technion.ac.il}
\urladdr{http://www.math.technion.ac.il/~pinsky/}

\subjclass[2000]{ 60J60} \keywords{}

\begin{abstract}
Consider a diffusion process corresponding to the operator $L=\frac12a\frac{d^2}{dx^2}+b\frac d{dx}$ and which is  transient
to $+\infty$. For $c>0$, we give an explicit criterion in terms of the coefficients $a$ and $b$ which determines whether or not
the diffusion almost surely eventually stops making down-crossings of length $c$. As a particular case, we show that
if $a=1$, then the diffusion almost surely stops making down-crossings of length $c$ if
$b(x)\ge\frac1{2c}\log x+\frac\gamma c\log\log x$, for some $\gamma>1$ and for large $x$, but makes down-crossings of length $c$ at arbitrarily large
times if $b(x)\le\frac1{2c}\log x+\frac1c\log\log x$, for large $x$.

\end{abstract}
\maketitle
\section{Introduction and Statement of Results}\label{Intro}

Consider the one-dimensional diffusion process $X(t)$ on $R$  corresponding
to the operator $L=\frac12a(x)\frac{d^2}{dx^2}+b(x)\frac{d}{dx}$, where $a$ is continuous and $b$ is locally bounded and measurable.
The question we address in this note is as follows: for $c>0$ and a given $a=a(x)$, how large a drift $b=b(x)$ is needed  in order for the process
to almost surely eventually stop making down-crossings of length $c$?

Let $P_x$ denote the measure
on $\Omega\equiv C([0,\infty),R)$, the space of continuous functions from $[0,\infty)$ to $R$, induced by
the above diffusion starting from $x\in R$.
Denote functions in $\Omega$ by $\omega=x(\cdot,\omega)$.
For each $c>0$, define the set-valued function  $S_c$ on $\Omega$ by
\begin{equation}
S_c(\omega)=\{x(t,\omega):t\ge0\ \text{and}\ \exists s>t\
\text{such that}\ x(s,\omega)\le x(t,\omega)-c\}.
\end{equation}
We refer to $S_c(\omega)$ as the  \it $c$-down-crossed range of $\omega$.\rm\
Let
$$
\sigma_c(\omega)=\inf\{t\ge0:\exists s>t\ \text{such that}\ x(s,\omega)\le x(t,\omega)-c\}
$$
 and let
\begin{equation*}\label{onset}
l_c(\omega)=\begin{cases}& x(\sigma_c(\omega),\omega),\ \text{if }\ \sigma_c(\omega)<\infty\\
&\infty,\ \text{if}\ \sigma_c(\omega)=\infty\end{cases}
\end{equation*}
denote the \it   onset location \rm\ of $S_c(\omega)$.
If the diffusion process $X(t)$ is recurrent, then trivially
$P_x(S_c=R)=1$, while if the diffusion almost surely converges to $-\infty$,
then clearly $P_x(S_c=(-\infty,\mu])=1$, where $\mu=\sup_{t\ge0}X(t)$.
For this reason, we will  be interested in the case
that the diffusion almost surely converges to $+\infty$.
As is well-known \cite{P}, this occurs if and only if
\begin{equation}\label{transient}
\int_{-\infty}\exp(-\int_0^x\frac {2b}a(y)dy)dx=\infty\ \
\text{and}\ \int^\infty\exp(-\int_0^x\frac {2b}a(y)dy)dx<\infty.
\end{equation}

Define
\begin{equation}\label{u}
 u(x)=\int_0^x\exp(-\int_0^y\frac{2b}a(z)dz)dy.
\end{equation}
Note that $u$ is harmonic for $L$; that is $Lu=0$.
\begin{theorem}\label{basic}
Let $P_x$ be the measure on $C([0,\infty),R)$ induced by  the diffusion process $X(t)$ starting from $x$ and
corresponding to the operator $L=\frac12a\frac{d^2}{dx^2}+b\frac d{dx}$.
Let $c>0$.

\noindent i.  The onset location $l_c$ of the $c$-down-crossed range $S_c$ satisfies
$$
P_x(l_c>x+\gamma)=\exp(-\int_x^{x+\gamma}\frac{u'(y)}{u(y)-u(y-c)}dy), \ \gamma>0,
$$
where $u$ is as in \eqref{u}.

\noindent ii. Assume that the diffusion is transient to $+\infty$; that is, assume that
\eqref{transient} holds. Then the $c$-down-crossed range $S_c$ is almost surely bounded, or equivalently, the diffusion $X(t)$
almost surely eventually stops
making $c$-down-crossings, if and only if
\begin{equation}\label{criterion}
 \int^\infty\frac{u'(x)}{u(x)-u(x-c)}dx<\infty.
\end{equation}

\end{theorem}

As an application of Theorem \ref{basic}, consider Brownian motion with a drift, corresponding to the operator
$L=\frac12\frac{d^2}{dx^2}+b\frac d{dx}$. Then $X(t)$ satisfies the stochastic differential equation
$X(t)=B(t)+\int_0^tb(X(s))ds$, where $B(t)$ is a Brownian motion. We have the following result.

\begin{theorem}\label{drift}
Consider the diffusion process $X(t)$ corresponding to
the operator $L=\frac12\frac{d^2}{dx^2}+b\frac d{dx}$.
Let $c>0$.

\noindent i. If $b(x)\ge\frac1{2c}\log x+\frac\gamma{ c}\log\log x$, for some $\gamma>1$
and for sufficiently large $x$, then the diffusion $X(t)$ a.s.
 eventually stops making $c$-down-crossings; that is, the $c$-down-crossed range $S_c$ is  bounded a.s.;

\noindent ii. If $b(x)\le \frac1{2c}\log x+\frac1{c}\log\log x$ for sufficiently large $x$, then the diffusion $X(t)$ a.s. makes
$c$-down-crossings for arbitrarily large $t$; that is, the $c$-down-crossed range $S_c$ is  unbounded a.s.;

\noindent iii. If $\liminf_{x\to\infty}\frac{b(x)}{\log x}=\infty$, then the diffusion $X(t)$ a.s. eventually stops making $c$-down-crossings for all $c>0$;
that is, the $c$-down-crossed range $S_c$ is bounded for all $c>0$ a.s.;

\noindent iv. If $\limsup_{x\to\infty}\frac{b(x)}{\log x}=0$, then the diffusion $X(t)$ a.s.  makes $c$-down-crossings for arbitrarily large $t$ for all $c>0$;
that is, the $c$-down-crossed range $S_c$ is unbounded for all $c>0$ a.s..

\end{theorem}
\noindent \bf Remark.\rm\
For a one-dimensional diffusion transient to $+\infty$, a general way to measure the tendency of the process to down-cross
is to ask for which unbounded, increasing sequences $\{d_m\}_{m=1}^\infty$ one has
 $P(A^{\{d_m\}}_n$ i.o.$)=1$ and for which such sequences one  has $P(A^{\{d_m\}}_n$ i.o.$)=0$,
where $A^{\{d_m\}}_n$ is the event that the process down-crosses from
$d_n$ to $d_{n-1}$.
It follows immediately from Theorem \ref{drift} that if $b(x)=\frac1{2c}\log x$, for large $x$, then when $d_m=\gamma m$, one has
$P(A_n^{\{d_m\}}$ i.o.$)=0$, if $\gamma>c$; furthermore as would be expected from the theorem, one can also prove that
$P(A_n^{\{d_m\}}$ i.o.$)=1$, if $\gamma\le c$.

It is interesting to contrast the above down-crossing behavior with the down-crossing behavior
in the case that the drift is of ``diffusion type,'' that is of the order $O(\frac1x)$ as $x\to\infty$.
Consider the Bessel process on $(0,\infty)$ corresponding to the operator $\frac12\frac{d^2}{dx^2}+\frac{k-1}{2x}\frac d{dx}$, and
assume that   $k>2$ so that the process is transient to $+\infty$. Of course, if $k$ is an integer, then the process
corresponds to the absolute value of $k$-dimensional Brownian motion.
For $\rho>0$, consider  $d_m=(m!)^\rho$ and denote $A_n^{\{d_m\}}$ by
  $A_n^\rho$.
The proof of Proposition 2 in  \cite{BP} shows that
$P_x(A_n^\rho \ \text{i.o.})=\begin{cases}& 0,\ \text{if} \ \rho>\frac1{k-2};\\
&1,\ \text{if}\ \rho\le\frac1{k-2}.\end{cases}$
\section{Proofs}

\noindent \bf Proof of Theorem \ref{basic}. \it (i)\rm\
Let $x\in R$ and let  $\gamma>0$. For $n$ a positive integer, define $x^{(n)}_k=x+\frac{k\gamma}n$, $k=0,1\cdots$.
Let $\tau_r=\inf\{t\ge0: X(t)=r\}$. By the strong Markov property, one has
\begin{equation}\label{inequality}
\prod_{k=0}^{n-1}P_{x^{(n)}_k}(\tau_{x^{(n)}_{k+1}}<\tau_{x^{(n)}_{k+1}-c})\le P_x(l_c>x+\gamma)\le
\prod_{k=0}^{n-1} P_{x^{(n)}_k}(\tau_{x^{(n)}_{k+1}}<\tau_{x^{(n)}_k-c}).
\end{equation}
As is well-known, since $u$ is harmonic, one has $P_w(\tau_z<\tau_y)=\frac{u(w)-u(y)}{u(z)-u(y)}$, for $y<w<z$.
Thus, since $u$ is uniformly Lipschitz on bounded intervals, one has
\begin{equation}\label{bigcalc}
\begin{aligned}
&\log\prod_{k=0}^{n-1} P_{x^{(n)}_k}(\tau_{x^{(n)}_{k+1}}<\tau_{x^{(n)}_k-c})=
\sum_{k=0}^{n-1}\log\frac{u(x^{(n)}_k)-u(x^{(n)}_k-c)}{u(x^{(n)}_{k+1})-u(x^{(n)}_k-c)}\\
&=\sum_{k=1}^{n-1}\log\left(1-\frac{u(x^{(n)}_{k+1})-u(x^{(n)}_k)}{u(x^{(n)}_{k+1})-u(x^{(n)}_k-c)}\right)=
-\sum_{k=1}^{n-1}\frac{u(x^{(n)}_{k+1})-u(x^{(n)}_k)}{u(x^{(n)}_{k+1})-u(x^{(n)}_k-c)}+O(\frac1n)\\
&=-\frac\gamma n\sum_{k=1}^{n-1}\frac{u'(z^{(n)}_k)}{u(x^{(n)}_{k+1})-u(x^{(n)}_k-c)}+O(\frac1n)
\ \text{as}\ n\to\infty, \ \text{where}\    x^{(n)}_k \le z^{(n)}_k\le x^{(n)}_{k+1}.
\end{aligned}
\end{equation}
Letting $n\to\infty$ in \eqref{bigcalc}, one obtains
\begin{equation}\label{limit}
\lim_{n\to\infty}\log\prod_{k=0}^{n-1} P_{x^{(n)}_k}(\tau_{x^{(n)}_{k+1}}<\tau_{x^{(n)}_k-c})=-\int_x^{x+\gamma}\frac{u'(y)}{u(y)-u(y-c)}dy.
\end{equation}
An almost identical calculation shows that the left hand expression in \eqref{inequality} also converges to the right hand side
of \eqref{limit} when $n\to\infty$. Thus, one concludes from \eqref{inequality} that
$P_x(l_c>x+\gamma)=\exp(-\int_x^{x+\gamma}\frac{u'(y)}{u(y)-u(y-c)}dy)$.
\medskip

\it\noindent (ii)\rm\ First assume that $\int^\infty\frac{u'(x)}{u(x)-u(x-a)}dx=\infty$.
Then by part (i) the onset location $l_c$ of $S_c$ is a.s. finite. Since the diffusion is transient to $+\infty$,
after it makes a $c$-down-crossing from the level $l_c$, it
will a.s. return to the level $l_c$. Starting anew from $l^{(1)}_c\equiv l_c$, by part (i) the diffusion will again a.s. make a $c$-down-crossing
with some onset location $l^{(2)}_c>l^{(1)}_c$. Continuing in this way, it follows that the set $S_c$ is a.s. unbounded.

Now assume that
$\int^\infty\frac{u'(x)}{u(x)-u(x-a)}dx<\infty$.
Under $P_x$, the probability of ever making a $c$-down-crossing is $q_x
\equiv1-\exp(-\int_x^\infty\frac{u'(y)}{u(y)-u(y-c)}dy)\in(0,1)$.
 If a down-crossing is made, with onset location $l_c\ge x$,
 then since the diffusion is transient
to $+\infty$, it will a.s.  eventually  return to $l_c$. Starting anew from $l_c$, the probability of making another $c$-down-crossing
is $1-\exp(-\int_{l_c}^\infty\frac{u'(y)}{u(y)-u(y-c)}dy)\le q_x$. Continuing like this, it follows that the diffusion will a.s. stop making $c$-down-crossings.
\hfill $\square$

\bigskip

For the proof of Theorem \ref{drift}, we will need a monotonicity result. If $L_i=\frac12a(x)\frac{d^2}{dx^2}+b_i(x)\frac{d}{dx}$,
for $i=1,2$, with $b_1\le b_2$, then a well-known coupling shows that the process corresponding to
$L_2$ and  starting at some $x$ stochastically dominates the process corresponding to $L_1$ and starting from the same point $x$.
It seems intuitive that in such a case, the number of $c$-down-crossings of the process corresponding to $L_1$ should stochastically
dominate the number of $c$-down-crossings of the process corresponding to $L_2$, however an appropriate, simple coupling doesn't
seem obvious. We obtain such a monotonicity result by Frech\'et differentiating the integrand in \eqref{criterion}.

\begin{proposition}\label{comp}
Let $L_i=\frac12a(x)\frac{d^2}{dx^2}+b_i(x)\frac{d}{dx}$, $i=1,2$, with $b_1\le b_2$. If the diffusion corresponding to
$L_1$ eventually stops $c$-down-crossing, then so does the diffusion corresponding to $L_2$.
\end{proposition}
\begin{proof}
Let $H(b,x)=\frac{u'(x)}{u(x)-u(x-c)}$, where $u$ is as in \eqref{u}.
By the assumption in the proposition, the diffusion corresponding to $L_1$ must be transient to $+\infty$. Since
$b_2\ge b_1$ the same holds for the diffusion corresponding to $b_2$.
From Theorem \ref{basic}, it suffices to show that $H(b_1,x)\ge H(b_2,x)$.
To show this, it suffices to show that if $q=q(x)\ge0$, then the Frech\'et derivative
$H_q(b,x)=\lim_{\epsilon\to0}\frac{H(b+\epsilon q,x)-H(b,x)}
{\epsilon}\ge0$.
One calculates that
\begin{equation*}
H_q(b,x)=\frac{\exp(-\int_0^x\frac{2b}a(y)dy)\int_{x-c}^x\exp(-\int_0^y\frac{2b}a(z)dz)(\int_y^x\frac{2q}a(z)dz)dy}
{(\int_{x-c}^x\exp(-\int_0^x\frac{2b}a(y)dy))^2}\ge0.
\end{equation*}
\end{proof}

\noindent\bf Proof of Theorem \ref{drift}.\rm\
Parts (iii) and (iv) follow immediately from parts (i) and (ii).
In light of Proposition \ref{comp}, to prove parts (i) and (ii), it suffices to consider the integral in \eqref{criterion} with
 $b(x)= \frac1{2c}\log x+\frac\gamma c\log\log x$, for large $x$, and show that this integral is infinite if
$\gamma=1$ and finite if $\gamma>1$. (For $b$ of this form \eqref{transient} holds.)

We  apply l'H\^opital's rule to the quotient
\begin{equation}\label{equal}
\frac{(x\log^{2\gamma-1} x) u'(x)}{u(x)-u(x-c)}=
\frac{(x\log^{2\gamma-1} x) \exp(-\int_0^x2b(y)dy)}{\int_{x-c}^x\exp(-\int_0^y2b(z)dz)dy}.
\end{equation}
It is clear that the numerator and denominator of the right hand side of \eqref{equal} tend to 0 as $x\to\infty$.
Differentiating and doing some algebra, we obtain
\begin{equation}\label{lop}
\frac{\left((x\log^{2\gamma-1} x)\exp(-\int_0^x2b(y)dy) \right)'}{(\int_{x-c}^x\exp(-\int_0^y2b(z)dz)dy)'}=
\frac{2b(x)(x\log^{2\gamma-1} x)+\text{lower order terms}}{\exp(\int_{x-c}^x2b(y)dy)-1}.
\end{equation}
Some standard analysis shows that
\begin{equation}\label{calculations}
\begin{aligned}
&\int_{x-c}^x \frac1 c\log ydy=\log x+2+o(1), \ \text{as}\ x\to\infty,\\
&\int_{x-c}^x\frac{2\gamma}c\log\log ydy=2\gamma\log\log x+o(1),\ \text{as}\ x\to\infty.
\end{aligned}
\end{equation}
Using \eqref{equal}-\eqref{calculations} along with the fact that $2b(x)=\frac1c\log x+\frac{2\gamma}c\log\log x$, for large $x$,
we obtain
\begin{equation*}
\lim_{x\to\infty}\frac{\left((x\log^{2\gamma-1} x) u'(x)\right)'}{\left(u(x)-u(x-c)\right)'}=\frac1{ce^2}.
\end{equation*}
It then follows from l'H\^opital's rule that
$\frac{ u'(x)}{u(x)-u(x-c)}\sim\frac1{ce^2}\frac1{x\log^{2\gamma-1}x},\ \text{as}\ x\to\infty$.
Thus, $\int^\infty\frac{u'(x)}{u(x)-u(x-c)}dx$ is infinite if $\gamma=1$ and finite if $\gamma>1$.
\hfill $\square$


\begin{thebibliography}{9}

\bibitem{BP} Ben-Ari, I. and  Pinsky, R. G.
\emph{ Absolute continuity/singularity and relative entropy properties for
probability measures induced by diffusions on infinite time intervals},
 Stochastic Process. Appl. \textbf{115} (2005), 179--206.

\bibitem{P} Pinsky, R. G.,
\emph{Positive Harmonic Functions and Diffusion},
    Cambridge Studies in Advanced Mathematics
    \textbf{45}, Cambridge University Press,
      (1995).
\end{thebibliography}
\end{document}